\newcommand{\CR}{\hbox{{$\cal R$}}}
\newcommand{\CE}{\hbox{{$\cal E$}}}
\newcommand{\CC}{\hbox{{$\cal C$}}}
\newcommand{\CS}{\hbox{{$\cal S$}}}

\newcommand{\Dsl}{{D\!\!\!\!/}}

\newcommand{\cg}{\mathfrak{g}}
\newcommand{\swap}{\leftrightarrow}
\newcommand{\R}{\mathbb{R}}
\newcommand{\C}{\mathbb{C}}
 
\newcommand{\eps}{{\epsilon}}

\newcommand{\extd}{{\rm d}}
\newcommand{\del}{\partial}
\newcommand{\isom}{{\cong}}
\newcommand{\Ad}{{\rm Ad}}
\newcommand{\tens}{\mathop{\otimes}}
\newcommand{\la}{{\triangleright}}
\newcommand{\End}{{\rm End}}
\newcommand{\id}{{\rm id}}
\newcommand{\<}{\langle}
\renewcommand{\>}{\rangle}

\newcommand{\eqn}[2]{\begin{equation}#2\label{#1}\end{equation}}

\newcommand{\lcross}{{>\!\!\!\triangleleft}}

\renewcommand{\o}{{{}_{\scriptscriptstyle(1)}}}
\renewcommand{\t}{{{}_{\scriptscriptstyle(2)}}}
\newcommand{\bz}{{{}^{\scriptscriptstyle(\bar 0)}}}
\newcommand{\bo}{{{}^{\scriptscriptstyle(\bar 1)}}}
\newcommand{\bt}{{{}^{\scriptscriptstyle(\bar 2)}}}
\newcommand{\uo}{{{}^{\scriptscriptstyle(1)}}}
\newcommand{\ut}{{{}^{\scriptscriptstyle(2)}}}

\documentclass[10pt]{article}
\usepackage{amssymb,amsmath,epsfig}
\newtheorem{lemma}{Lemma}

\newtheorem{theorem}[lemma]{Theorem}

\newtheorem{corol}[lemma]{Corollary}

\begin{document}\baselineskip 14pt

{\ }\qquad \hskip 4.3in \vspace{.2in}

\begin{center} {\LARGE Conceptual Issues for Noncommutative Gravity on
Algebras and Finite Sets}
\\ \baselineskip 13pt{\ }\\
{\ }\\ Shahn Majid \\ {\ }\\ School of Mathematical Sciences,
Queen Mary and Westfield College\\ University of London, Mile End
Rd, London E1 4NS\footnote{Reader and Royal Society University
Research Fellow}

\end{center}

\begin{quote}\baselineskip 15pt
\noindent{\bf Abstract} We discuss some of the issues to be
addressed in arriving at a definitive noncommutative Riemannian
geometry that generalises conventional geometry both to the
quantum domain and to the discrete domain. This also provides an
introduction to our 1997 formulation based on quantum group frame
bundles. We outline now the local formulae with general
differential calculus both on the base `quantum manifold' and on
the structure group  Gauge transforms with nonuniversal calculi,
Dirac operator, Levi-Civita condition, Ricci tensor and other topics are also
covered. As an application we outline an intrinsic or relative theory of
quantum measurement and propose it as a possible framework to explore the link
between gravity in quantum systems and entropy. \end{quote}

\section{Introduction}    There has been a
lot of interest in recent years in a theory of  `Riemannian geometry' and
gravity applicable to general possibly  noncommutative, possibly discrete
algebras. The motivations are  obvious enough and have been featured in many a
research proposal  in the last several years. In this note I want to discuss
some of  the issues in arriving at the proper formulation as well as to
provide an introduction to my own solution to this problem some
years ago, and which I presented in detail at the Euroconference.
Probably the biggest problem with this formulation is that it is
technically too hard for most readers; I will try to address some
of that here. I also want to put in writing the announcements
made in my lectures of some of my results, which will appear in
the long paper \cite{Ma:rieq} in preparation.

We start with some motivation for such a theory, which in my case
(for about 15 years) has always been quantum gravity or Planck
scale physics. The first paper that I am aware of in which actual
noncommutative geometries based on quantum group methods were
proposed as a model of Planck scale physics was \cite{Ma:pla}. We
put forward here (in 1987) the view that:

\begin{quote}
A noncommutative noncocommutative Hopf algebra should be viewed
as a toy model of a system in which both quantum effects (the
noncommutativity) and gravitational curvature effects (the
noncocommutativity) are unified.
\end{quote}

Many years later, with a theory of quantum Riemannian geometry,
we can now actually measure the curvature in a (possibly) quantum
system, for example relating it to the noncocommutativity in the
case of a Hopf algebra or quantum group, and thereby make this
idea precise. Many inspired works were born from  semantic confusions
and here I want to note for example that if one views one of
Drinfeld's quasitriangular Hopf algebras $U_q(\cg)$ this way (as
a `coordinate algebra', which is not the usual way) then one  can
expect the following identity
\eqn{RieDri}{ R_{\rm Riemann}\sim \CR_{\rm Drinfeld},}
where the latter is the quasitriangular structure or universal
R-matrix that controls the noncocommutativity. The theory to be
described below has the power to make this precise.

More physically, the paper \cite{Ma:pla} introduced some concrete
ideas for the Planck scale relating to a duality between quantum
theory and geometry which could be viewed as a modern version of
Mach's principle, and such ideas could be explored now in more
detail. One should certainly like to have more realistic and
relativistic models with similar phenomena. Also in general
terms, there are many indications of a profound relationship
between gravity, statistical mechanics and the `problem of time'.
The tools of noncommutative Riemannian geometry would provide the
power to resolve those relationships in detail. We would be able
to deal with algebras which are actual quantum systems, compute
entropy of states etc, and at the same time compute their
geometrical structure, curvature etc, without classical limits,
and hence sort out these issues algebraically (I would guess that
they relate to some form of the above -mentioned duality). Again,
this has been explored in the very simple case of Hopf algebras,
at least a little, via `quantum random walks' on
them\cite{Ma:wal} but before the tools of Riemannian geometry
were available.

Another impetus is that string theorists are finally beginning to
accept this view (which was not true 10 years ago for example)
although the noncommutative tori found there (at the algebraic
level at which physicists tend to work) are too trivial to be
good candidates for quantum spacetime, being the usual Heisenberg
or Weyl algebra in one form or another.

In a different direction, one can use noncommutative geometry
methods to regularise even ordinary quantum field theories on
flat space. The first result of this kind I believe appeared in
1990 in \cite{Ma:reg} where it was proposed to consider the
enveloping algebra $U_q(b_+)$ `up side down' as a noncommutative
spacetime and shown that this leads to better convergence of some
integrals. In the limit $q=1$ our proposal was for $[x,t]=x$ to
be spacetime. Recently such noncommutative `kappa-Minkowski'
spacetimes have been studied using the integral  on this algebra
in \cite{Ma:reg} and ensuing Fourier transform to make the
first mathematically meaningful predictions for gamma ray
propagation that could be tested in principle by cosmological
data\cite{AmeMa:wav}. Also many years later, regularisation by
q-deformation is  finally becoming a reality for viable quantum
field theories, as shown by Robert Oeckl's talk at this
conference; see his seminal paper\cite{Oec:bra}.

The improved properties for regularisation can be traced to
finite differences rather than usual differentiation. The natural
emergence of these in noncommutative geometry obviously makes any
successful Riemannian geometry of that type also a natural and
less ad-hoc replacement for lattice and other discrete modeling
methods for quantum gravity. One can either work with a
finite-difference calculus on classical (not lattice) spacetime
by choosing an appropriate noncommutative differential calculus,
or one can go further and make the space or spacetime itself
discrete or finite. This is another obvious application of a
successful noncommutative Riemannian geometry.

Also, Connes some years ago gave an approach to the standard
model of elementary particle physics in which conventional
internal gauge symmetries are recovered (along with the Higgs
potential) in a much more elegant manner than usually presented, by
extending spacetime by a finite-dimensional noncommutative
algebra or `discrete' part. The problem with this approach is
that some of the most interesting physical data relating to
masses in the standard model is encoded in the `Dirac operator'
on the finite algebra, but this can be practically any
self-adjoint operator. If we had a true quantum Riemannian
geometry that  included discrete or finite systems then we would
know what was the natural or canonical Dirac operator in a given
context from a geometrical point of view and thereby obtain
direct physical predictions about the standard model. So the same
technology would have a second quite different application to
physics. One can go further and relate it to quantum gravity as
well through Connes more recent `spectral action principle'.

If this list of motivations is not enough, I want to present a
final one. This is that one can gain geometric tools for ordinary
quantum mechanics and address fundamental interpretational
problems there. For what is the correspondence principle other
than the idea that certain operators in the quantum theory like
$x,p$ behave similarly enough to their classical counterparts
that we correspond them? Noncommutative geometry for the first
time offers a way to turn that very woolly idea into something
precise, i.e. to extend the correspondence principle itself into
something rather more sophisticated in which we say that a given
quantum algebra has this or that quantum geometric structure and
therefore identify correctly the pieces for its classical
geometric limit. Actually this is not something just for
philosophers of quantum mechanics; if we are going to achieve and
understand quantum gravity we will probably need to fully
understand the measurement process in a more intrinsic manner
than `somebody outside' making a measurement. This is somehow
tied up with the entropy questions we began the discussion with.

Let us now summarise some of the  requirements stemming from the
above motivations. I would identify 4 key demands for
noncommutative Riemannian geometry.
\begin{enumerate}
\item The theory should have full conceptual continuity with
conventional Riemannian geometry so that this is included as a
special case (along with finite geometry and noncommutative
algebras as other specialisations).
\item There should be either a global picture (for nontrivial topology)
of the whole theory or local gauge transformations from which a global
picture could be patched.
\item The theory should be powerful enough to include a wide variety
of genuine examples covering a full range of objects that people
would want to consider `geometrically' and with similar or more
degree of `flabbiness' as usual geometry.
\item The theory should  probably in principle include spinors since
these are needed for matter fields.
\end{enumerate}

To give an example of a theory failing 2., supermanifolds are
basically ordinary manifolds with a `Grassmann' part glued on;
the latter is always flat and not really capable of its own
topological structure. For an example of an approach failing 3.,
the `idea' that vector fields should be derivations is obvious
(and makes sense for any noncommutative algebra) but works only
for very special matrix type algebras; in many other cases, (e.g.
for quantum groups), the natural partial derivatives are not
derivations at all, and indeed essentially none exist.  This
problem may also occur in Connes' spectral triple approach where,
without a full range of noncommutative examples, one cannot be
sure about the applicability of the axioms even if they are
correct in the commutative case and mathematically elegant. That
is why I want to emphasise 3.

Here at this conference several interesting talks were related to
these topics. Sorace made some steps towards Dirac operators on
Euclidean quantum groups, while Castellani made an attempt at a
theory of gravity on finite groups. None of the works, even in the finite
or discrete case, however,  fit  our primary demand of conceptual
containment and continuity with classical geometry, either missing key
definitions such as compatibility between the so-called `spin connection' and
the metric, or using ad-hoc definitions heavily dependent on a lattice or
group structure on the spacetime, or limited to some form of universal calculus
much bigger than classical. Also on point 4. it should also be said that
essentially all  attempts at
noncommutative Riemannian geometry theories other than the frame  bundle
approach introduced in \cite{Ma:rie} have been based on  positing  the
covariant derivative $\nabla$ directly via certain  axioms, i.e. a theory of
linear connections. Such an approach can  never interact properly with matter
fields, particularly spinors,  since these need at their base a principal
bundle of some sort to  which vector bundles such as the spin bundle are
associated. The  only way that has been found to circumvent that is Connes's
approach to jump over covariant derivatives altogether and posit  the `Dirac
operator` directly as the basic object, as mentioned  above. However, in this
approach it is hard to pick out  connections or much else of the
infrastructure of differential  geometry at all, other than in the commutative
case, i.e. one  cannot test our primary `continuity' requirement around the
commutative case.  Connes' axioms may in fact have to be  significantly
modified to reach the standard q-deformation  examples, for instance.

I would like to stick my neck out and claim that our 1997 frame  bundle
approach to Riemannian geometry in \cite{Ma:rie} on the  other hand {\em does}
already more or less meet all four criteria  at least in principle. The paper
itself focussed on two extremes,  checking continuity with the classical case
and at the other  extreme the universal differential calculus, but it also
indicated the theory for general calculi using the principal  bundle theory
for general differential calculi already known by  then in \cite{BrzMa:gau},
i.e. this was a matter of presentation  and not a fundamental or conceptual
problem with the theory.  Basically, the examples had not been worked out for
nonuniversal  calculi to justify a full presentation, an oversight that will
be  filled in the forthcoming paper\cite{Ma:rieq} with examples for  all
standard quantum groups and their standard calculi, as  announced here in
Torino.

It should also be clear in all this the role of quantum groups
themselves: I would view them as nothing but a good playground of  `naturally
arising' and evidently geometrical objects, just as  Lie groups provided the
perfect playground for the development of  modern classical differential
geometry in the late 19th and early  20th century. Particularly, q-deformed
examples, because of their  parameter, allow us to test the continuity of our
noncommutative  Riemannian geometry on an `open set of models' near the
commutative point. This goes much beyond simply having the right  answer in
the commutative case. However, any theory of  `Riemannian geometry' that works
{\em only} for quantum groups,  or only for finite sets, etc., will not meet
our key demand and  will most likely be viewed as ad-hoc when the full theory
is  known. This excludes all sorts of rigid theories based on  integrable
systems etc; they should be examples of some general  theory but are not that
theory.

In fact, some examples that one would wish to include already
force one to generalise away from a quantum group as `gauge
group' to a general algebra (or rather, coalgebra in our dual
language) in its role. It turns out that most of what we want to
do can be done (with a lot more care) at that much more general
level\cite{BrzMa:geo} but along exactly the same lines as in
\cite{Ma:rie}. This step is certainly important to meet our
demand 3. and we will say more about it in the last section, but
at the moment only the special case of a quantum group framing
(which could of course be a classical or finite group if one
wants) has been more or less fully worked out with general
differential calculus.

\section{Size matters for differential calculi}

We set the scene with a brief discussion of the choice of
differential calculus, which in common concensus means a choice
of exterior algebra $\Omega(M)$ for our `coordinate' algebra $M$.
$M$ could be the functions on a finite group, a matrix algebra, a
quantum group etc. As in conventional geometry one can focus on
each degree of differential forms one at a time, i.e. first
define $\Omega^1(M)$ and let the higher order be unspecified
until one needs them. More formally, one can say that the higher
order just have the relations among higher order forms implied by
Leibniz and $\extd^2=0$ from what was assumed at degree 1. This
is the `maximal prolongation'. On the other hand, that is not at
all what would give the right answer in the classical case. There
will typically be an infinite exterior algebra with forms of all
degrees and no top form, etc., because sufficiently many
relations at higher degree are not demanded by Leibniz and
$\extd^2=0$ alone.

At the moment there are only three ways to specify a reasonable calculus to
all orders of differentials. In either case one starts with the
universal exterior algebra $\Omega M$ associated canonically to
any algebra -- this has very little structure at all and is far
far too big to be classical. One then specifies the quotients at
each degree by some scheme to get down to a calculus of
`reasonable' size at each order. For example, at degree one, the
universal calculus $\Omega^1 M$ is the subspace of $M\tens M$
given by the kernel of the product map. Its universal $\extd$ is
\eqn{univd}{ \extd m=1\tens m-m\tens 1,\quad \extd: M\to \Omega^1
M.} It is a bimodule so we can multiply forms by `functions' $M$
from the left or right. Any other choice $\Omega^1(M)$ is
specified by quotienting, i.e. by setting to zero some
subbimodule $N_M\subset \Omega^1 M$. In almost all cases, e.g
even for finite sets, the quotients will remain bimodules with
different left and right multiplications, i.e. `quantum' in the
sense that forms and functions do not commute. In this sense the
commutative calculus that we are familiar with classically is
very artificial and atypical.

The problem is that one needs a scheme to specify this
quotienting to all orders or at least to the orders needed for
physics (which for gauge theory means mainly to order 2 or 3).
One scheme is to mimic classical ideas. For example one can have
a `discrete geometry' (where $M$ is the algebra of functions on a
discrete set) by starting with some existing manifold or more
generally a sigma-algebra and use the intersection structure of a
good open cover to specify a calculus on the indexing set of the
cover, see for example\cite{BrzMa:dif} (which differed from other
similar proposals). This would give an example of a discrete
version (using noncommutative geometry for the differentials) of
your favorite classical manifold.

Another scheme is to suppose that the space has a group
structure, i.e. $M$ is a quantum group (or classical group
coordinate algebra). The coproduct $\Delta:M\to M\tens M$
expresses the group structure of course. For on a Lie group there
is a unique differential structure that is translation covariant
from the left and right and it specifies forms to all orders.
This was formulated many years ago by Woronowicz\cite{Wor:dif}
who showed in particular that once $\Omega^1(M)$ was fixed in a
bicovariant manner then there was a natural $\Omega(M)$ of the
right kind of size (reducing correctly in the Lie group case).
The basic idea is to set to zero those products of forms that
would be invariant under a generalised transposition or
`braiding' operator $\Psi$. The braiding is that of the category
of representations of the Drinfeld quantum double $D(M)$.
Recently by results in \cite{Ma:cla}\cite{BegMa:dif} and a
twisting theorem in \cite{MaOec:twi} the classification problem
for the remaining freedom of choice of $\Omega^1(M)$ for all main
classes of quantum groups was solved. So by now we know the menu
of choices of differential calculi in the quantum group case more
or less completely.

The third scheme is through Connes spectral triple, i.e. to
choose a more or less arbitrary operator $\Dsl$ that you would
like to be able to call Dirac acting on a vector space of the
form $M\tens W$ ($W$ would be the spinor space) and let it
determine the smallest differential calculus on $M$ that would be
needed for this to be the case. Most of Connes axioms involving
Hilbert spaces etc are not needed for this and one can formulate
it all as a purely algebraic construction. It is analysed a
little in comparison with the Woronowicz approach in
\cite{Ma:rieq}. Basically, consider $\pi:\Omega M\to \End(M\tens
W)$ defined for example on 1-forms by $\pi(m\tens n)=m[\Dsl,n]$
where $m,n$ act by multiplication. Then\cite{Con:geo} quotient
the universal $\Omega^1 M$ by $\ker\pi$ to define $\Omega^1(M)$.
Similarly at higher order we quotient by the differential ideal
generated by the kernel of $\pi$. This gives the right answer
classically if one really starts with the Dirac operator but all
depends on the choice of $\Dsl$. For a given algebra you may have
little or no idea which operator to take to get something finite
or resembling classical geometry.

To give an example, for polynomials in one variable over a field
$k$ the coirreducible translation invariant calculi have the form
\[ \Omega^1=k_\lambda[x],\quad  \extd
f(x)=\frac{f(x+\lambda)-f(x)}{\lambda},\quad f(x)\cdot
\alpha(\lambda,x)=f(x+\lambda)\alpha(\lambda,x),\quad \alpha(\lambda,x)\cdot
f(x)=\alpha(\lambda,x)f(x)\] for functions $f$ and one-forms $\alpha$.
Here $k_\lambda$ is a field extension of the form $k[\lambda]$
modulo $m(\lambda)=0$ and $m$ is an irreducible monic
polynomial.  For example, the most important field extension in
physics, $\R\subset\C$, can be viewed
noncommutative-geometrically with complex functions $\C[x]$ the
quantum 1-forms on the algebra of real functions
$\R[x]$\cite{Ma:fie}. There is nontrivial quantum DeRahm
cohomology in this case.

\section{To gauge or not to gauge; the global approach}

Once you have a differential calculus you can do cohomology, of
course. You can also do what I like to call $U(0)$ gauge theory
in which we regard a connection or gauge field simply as a 1-form
$A$ but transforming as \eqn{U(0)}{ A^\gamma=\gamma^{-1}A\gamma
+\gamma^{-1}\extd \gamma} under `gauge transform' by $\gamma$ an
invertible element of $M$. Or a unitary element in a
representation of the algebra, say. From a noncommutative
geometric point of view in which $M$ is `functions' this is a
trivial gauge theory even though the curvature \eqn{F(0)}{
F=\extd A+A\wedge A} transforms nontrivially by conjugation.

To do gravity on $M$ we need a full nonAbelian gauge theory in
which the gauge field has values in a general quantum (or
classical) group $H$. Such a theory was introduced in 1992 in
\cite{BrzMa:gau} and remains the only properly formulated global
and nonAbelian gauge theory for a quantum space that I am aware
of. More precisely we equip $H$ with a differential calculus
$\Omega^1(H)$ and consider the space $\Omega_0$ of left-invariant
1-forms as the dual of its Lie algebra in some sense, and gauge
fields should be maps from $\Omega_0$ (i.e. have values in
$\Omega_0^*$ in more conventional terms). We also need an algebra
$P$ with its own but compatible differential structure
$\Omega^1(P)$ which plays the role of the total space of a
principal bundle and on which $H$ (co)acts by a coaction
$\Delta_R$ with fixed subalgebra $M$ and some axiom of `local
triviality'. Remarkably, it is possible to do away with all of the
usual mess with local trivialisations and local charts (which we
will not have for a general algebra) and replace them by an
algebraic condition, the exact sequence
\eqn{exactness}{ 0\to
P\Omega^1(M)P\to\Omega^1(P){\buildrel{\rm ver}
\over\longrightarrow} P\tens \Omega_0\to 0.}
Here ${\rm
ver}=(\cdot\tens\id)\Delta_R$ defines for each element of
$\Omega_0^*$ the corresponding vertical vector field
$\Omega^1(P)\to P$ along the fiber and the exact sequence says
that the forms that they kill should be exactly the ones pulled
back from the base. A connection in the bundle is a splitting of
$\Omega^1(P)$ and can be cast as a connection form \eqn{connP}{
\omega:\Omega_0\to \Omega^1(P)} with properties just like in any
good book on differential geometry.

So the global theory exists, a nontrivial example, the
$q$-monopole, also in \cite{BrzMa:gau} showed that it worked, and
some generalisations for other example have been made also. But
are we happy after 8 years? Not completely because in physics we
also think we want to see local gauge transformation formulae
like in (\ref{U(0)}) at least for trivial bundles. And here there
is a problem. The problem is that the abstract mathematical and
global way works  but the local formulae that you know and love
have a lot more classical assumptions hidden in them than meet
the eye and which have to be appropriately generalised in a
smooth manner. Basically, let us consider what {\em should} be
the nonAbelian analogue of (\ref{U(0)}) which is a map
\eqn{connA}{ A:\Omega_0\to \Omega^1(M)} We can consider gauge
transform by $\gamma:H\to M$ in a similar manner but we wont be
able to `conjugate' $\Omega_0$ in the correct way. For to make
(\ref{U(0)}) we would want to apply the coproduct twice as a map
$\Omega_0\to H\tens\Omega_0\tens H$, keep the middle output for
$A$ and the outer two for $\gamma^{-1}$ and $\gamma$. This is an
analogue of conjugating the Lie algebra by the group, but it is
not one that works for a Hopf algebra. There is no such map
unless we take the universal differential calculus. If we look at
what happens classically what we really need is the adjoint
action of the group on the Lie algebra. Here we are in luck,
there is a well-defined adjoint coaction $\Ad:\Omega_0\to
\Omega_0\tens H$. We use it, apply $A$ and $\gamma$ and the
bimodule structure, \eqn{Aad}{ A^\gamma=\cdot\circ(A\tens
\gamma)\circ\Ad.}
The problem is that then the curvature is not
invariant as soon as the differential calculus is noncommutative
because the $\gamma^{-1}$ buried inside $\gamma\circ\Ad$ is on
the wrong side of $A$ in $\Omega^1(M)$. This is the fundamental
obstruction that has been around for many years now.

Yet there is no problem in the nonAbelian case with the full
global picture! After thinking about this for some years the
resolution of the problem is as follows\cite{Ma:rieq}. Firstly,
we can postpone the problem. After all, what do we want local
gauge transformations for anyway except to make sure things patch
up to a global picture? Having a global picture we {\em are}
already free to make bundle automorphisms or gauge
transformations as transformations on $P$ -- the problem is just
the local picture not the concept of gauge transformation itself.
It is only in the local picture that formulae like (\ref{U(0)})
come up.  Next, as we elaborate more of the local theory we will
see the gauge transformations appearing and eventually get a
local formulae like (\ref{U(0)}), which must exist because global
bundle automorphisms exist. To give a rough idea what the
eventual local formulae will look like (and to see why the result
will be rather too complicated in general just to guess) suppose that you have
a single global patch described by a trivial bundle. Such bundles
were considered in \cite{BrzMa:gau} and one example, dependent
already on the choice of trivialisation or `gauge' is the tensor
product bundle $P=M\tens H$. As already understood there, the new
feature in quantum groups gauge theory is that under a gauge
transformation $\gamma$ there will be a new trivial bundle with
the same vector space $M\tens H$ but in fact $P$ will no longer
be a tensor product but a cocycle cross product
\eqn{Pgamma}{P^\gamma=M_{\chi_\gamma}\lcross H.} So the very
algebraic form of the bundle changes under a gauge
transformation. The cross product and cocycle would all be
trivial if $M$ were commutative so you would never leave the class
$M\tens H$ in the classical case; the cocycle or `anomaly' is a
purely quantum effect. The new bundle is equivalent but our naive
description as a tensor product is not gauge invariant! So there
will be a new local picture for $A^\gamma$ and a familiar gauge
theory but only when we consider the whole class of cocycle
trivial bundles $P=M_{\chi}\lcross H$ stable under such $\gamma$
(with $\chi$ also transforming). There is no problem of principle
here at all except that the result is more complicated than we
would have liked. With this in mind, we will shortly give all the
formulae of quantum Riemannian geometry in `tensor product gauge'
but one should not try to make gauge transformations of them
without being prepared to do it in the manner above with a
cocycle. Also note that not all trivial bundles are gauge
equivalent to a tensor product one (it depends on the vanishing of
a nonAbelian cohomology where $\chi$ lives, see \cite[Ch. 6]
{Ma:book}), so we will be seeing in this way only the
cohomologically trivial sector of the global theory (other sectors
never visible classically are also possible).

In effect, our interim solution to the gauge problem is to work
globally or work in what I propose to call `tensor product
gauge'. Because the global theory exists  it is enough to work
(for trivial bundles) in one gauge like this, and consider other
algebraic gauges later. It still takes a lot of
work\cite{Ma:rieq} to prove that the local formulae indeed lead
to a global bundle on $P=M\tens H$ -- we are simply validating
them in a different way than by gauge transformation.

And if one really wants gauge transformations at all cost then
one can of course use the universal differential calculus on $H$.
Then $\Omega_0\subset H$ (it is the kernel of the counit) and one
can use the coproduct and a projection to define $\Omega_0\to
H\tens\Omega_0\tens H$ as we had first wanted to do. This was
also introduced in 1992 in \cite{BrzMa:gau} where all usual gauge
formulae are given (and this works for general $\Omega^1(M)$).
But such a theory in more usual terms would be allowing the gauge
field and curvature to have values in the whole quantum group
dual to $H$ as explained in \cite{BrzMa:gau}. Such a theory is OK
but it violates our principal demand of continuity with classical
geometry wherin $A$ and $F$ should have values in some kind of `Lie
algebra' and not in the enveloping algebra to resemble classical
geometry. It is only to achieve this that we have to work much
harder as explained above.

\section{Quantum Riemannian geometry -- global version}

Before passing to tensor product gauge, let us complete the
formalism of global quantum Riemannian geometry as basically
introduced in \cite{Ma:rie}. Firstly, in the theory of quantum
principal bundles one has the `associated bundles'
$\CE^*=\hom^H(V,P)$ and $\CE=(P\tens V)^H$ (the invariant
subalgebra of $P\tens H$) defined by any space $V$ in which $H$
(co)acts. Any connection on $P$ then induces covariant
derivatives on these, \eqn{covD}{ D:\CE\to \Omega^1(M)\tens_M\CE}
etc. The new ingredient for Riemannian geometry is a soldering
form $\theta:V\to P\Omega^1(M)$ such that the induced map $\CE\to
\Omega^1(M)$ by applying $\theta$ and then the bimodule structure
of $\Omega^1(P)$ is an isomorphism. What this does is to express
the cotangent bundle as associated to the principal bundle $P$.
This in my view captures the essence of what a manifold {\em is}
since all the ideas of charts etc are taken care of by the local
triviality of the bundle, which in turn we have expressed
globally by the exact sequence above.

We call a `quantum manifold' an algebra framed in this way. A
connection on $P$ then induces a covariant derivative
\eqn{nabla}{\nabla:\Omega^1(M)\to \Omega^1(M)\tens_M\Omega^1(M),}
which is just $D$ viewed under the framing isomorphism. A metric
for us is an element  \eqn{g}{ g\in
\Omega^1(M)\tens_M\Omega^1(M)} and corresponds to an isomorphism
$\CE^*\isom\CE$. This can be expressed in turn as the idea that
there is a second framing $\theta^*:V^*\to \Omega^1(M)P$ with
$V^*$ in the role of $V$, what I call the `coframing'. The
corresponding metric is $\theta^*\tens_P\theta$ applied to the
canonical element of $V^*\tens V$. So a `quantum Riemannian
manifold' is an algebra $M$ with a framing and coframing.

Finally, we should have compatibility conditions between the
framing and the connection. One is zero torsion and corresponds
to $\bar D\wedge\theta=0$ where $\bar D$ is a suitable left-handed
covariant derivative. The other is a `metric compatibility'
condition. This is the hardest to formulate and the solution
proposed in \cite{Ma:rie} required a slight generalisation of
Riemannian geometry itself in a manner appropriate to
non-symmetric metrics (we have not demanded symmetry of $g$
above) as `skew metric compatibility' or vanishing of `cotorsion'
$D\wedge\theta^*=0$.

The correspondence of this theory with conventional objects
(which we need for continuity with classical geometry)  takes the
form \eqn{corresp}{ g\swap \theta,\theta^*,\quad R\swap F,\quad
T\swap \bar D\wedge\theta,\quad \Gamma\swap D\wedge\theta^*} where
\eqn{curvR}{
R=((\id\wedge\nabla)-(\extd\tens\id))\circ\nabla,\quad
\Omega^1(M)\to \Omega^2(M)\wedge\Omega^1(M)} is the Riemann
curvature, \eqn{torT}{ T=\extd-\nabla\wedge,\quad
T:\Omega^1(M)\to \Omega^2(M)} is the torsion tensor and
\eqn{cotor}{  \Gamma=(\nabla\wedge\id
-\id\wedge\nabla)g+(T\tens\id)g,\quad \Gamma\in
\Omega^2(M)\tens_M\Omega^1(M)} is our novel concept of `cotorsion
tensor'. The difference between torsion and cotorsion evidently
expresses metric compatibility in a skew manner.

If one wants to be more conventional there is no problem writing
axioms for symmetry of the metric in a suitable form, such as
\eqn{symg}{\wedge(g)=0.} Also one could be more conventional and
ask for full metric compatibility under $\nabla$ in a more usual
way, however this does not seem to me very natural compared to the
elegance of the cotorsion approach in the frame bundle context
(where we keep everything symmetric between framing and
coframing) and moreover would still not achieve a unique
Levi-Civita connection; that uniqueness is a very special feature
of the frame bundle group classically being $O_n$ and cannot be
expected in general.

Now we can summarise and evaluate the degree of continuity of our
theory with conventional Riemannian geometry. This was a large
part of the work in \cite{Ma:rie}. The summary is that the
quantum Riemannian geometry that we propose reduces in the
classical case to a slight generalisation of conventional
geometry in which:

\begin{enumerate}

\item We allow any group $G$ in the `frame bundle', hence
a more general concept of a `frame resolution'
$(P,G,V,\theta_\mu^i)$ or {generalised manifold}. The choice of
group determines the range of $\nabla$ that can be induced.

\item The {generalised metric} $g_{\mu\nu}
=\theta^*_\mu{}^i\theta_{\nu i}$ corresponding to a coframing
$\theta^*_{\mu}{}^i$ is nondegenerate but need not be symmetric.

\item The {generalised Levi-Civita} connection defined as having
vanishing torsion and vanishing cotorsion respects the metric only
in a skew sense \eqn{cotorcla}{ \nabla_\mu g_{\nu\rho}-\nabla_\nu
g_{\mu\rho}=0} and need not be uniquely determined.

\end{enumerate}

We may still impose symmetry of the metric in some natural form
as above, etc., for strict continuity with classical Riemannian
geometry. In this instance, however, I think it is also a good
opportunity -- suggested by the noncommutative Riemannian
geometry -- to enlarge Riemannian geometry itself a little and
unify it with other branches of classical geometry. Thus, in the
other extreme when the generalised metric is totally
antisymmetric the vanishing of cotorsion minus torsion implies $\extd
g=0$ so $g$ is symplectic, as shown in \cite{Ma:rie}. This is
fully in keeping with what me might expect from T-duality and
other physical considerations as the first Planckian corrections
to classical geometry.

We also promised a good supply of examples. Announced in Torino
and to appear in \cite{Ma:rieq} we have

\begin{theorem}\cite{Ma:rie}\cite{Ma:rieq}
All quantum groups $M$ equipped with bicovariant calculi are
quantum manifolds in the above sense with framing by $H=M$ itself
and $\theta$ induced by the Maurer-Cartan form in \cite{Wor:dif}.
Moreover, all standard quantum groups $\C_q[G]$ with their
standard bicovariant differential calculi are quantum Riemannian
manifolds in the above sense with metric induced by the
braided-Killing form\cite{Ma:lie} on the braided-Lie algebra
associated to the differential calculus. \end{theorem}

From this we obtain in principle similar results for quantum
homogeneous spaces including spheres, planes etc. In fact, there is
a notion of comeasuring or quantum automorphism
bialgebra\cite{Ma:dif} for practically any algebra $M$ and when
this has an antipode (which typically requires some form of
completion) one can write $M$ as a quantum homogeneous space. So
almost any algebra $M$ is more or less a quantum manifold for
some principal bundle (at least rather formally and so far with
the universal calculus). When $M$ is equipped with more structure
to define a nonuniversal differential calculus in a systematic
way then this observation should extend to that level too, which
would then be  analogous to the idea that any classical manifold
is, rather formally, a homogeneous space of diffeomorphisms
modulo diffeomorphisms fixing a base point.

Also, let us mention the obvious global formulation of the Dirac
operator in the frame bundle approach. Given any other vector
space $W$ on which $H$ (co)acts we similarly have an associated
bundle $\CS$ say as explained above. The connection on the
principal frame bundle that induced the covariant derivative
$\nabla$ also induces a covariant derivative $D:\CS\to
\Omega^1(M)\tens_M \CS$. So the missing ingredient is just a
suitable map $\gamma:\Omega^1(M)\to \End(\CS)$ which locally
would be provided by $\gamma$ matrices. I am not confident about
the global formulation of this map, however the paper
\cite{Ma:rieq} contains examples, for instance on finite groups
using this bundle formulation and along these lines. i.e. the
theory includes spinors in principle and has been tested on some
examples.

Finally, we discuss the Ricci tensor and actions.  To make the
contraction for Ricci we can assume a `lifting map' or
bimodule inclusion \eqn{lift}{ i:\Omega^1(M)\tens_M\Omega^1(M)\to \Omega^2(M)}
so as to turn $R$ as a  2-form (with values in operators on 1-forms) into an
element of  $\Omega^1(M)\tens_M\Omega^1(M)$ (with values as before). After
that we can then take a suitable trace of $i(R)$ as a left $M$-module
endomorphism with values in the remaining (rightmost) copies of
$\Omega^1(M)\tens_M\Omega^1(M)$, which defines the Ricci tensor. Such a lift
$i$ is also the ingredient needed for an interior product on the exterior
algebra, so specifying this or a Hodge $*$ operation would also do the job in
a natural manner. How  these or $i$ are chosen depends on how the  higher order
differential calculus (or at least $\Omega^2(M)$) is defined but in the main
cases one has some natural proposals. Assuming $i$ has been chosen
we can then write down field equations via Ricci.  Contracting further,
we would obtain the scalar curvature. We  would still need an integral for an
action. Alternatively we could  define the action via $\Dsl^2$ in view of the
well-known  Lichnerowicz formula. In the finite case we could take a trace
while in general we could use the spectral action of Connes.  These are some
of the steps for which ambiguities exist in the  frame bundle approach. I
think there is no fundamental  obstruction, however, to completing the
programme. This is also  clearly reaching the point where it {\em should} tie
up with  Connes' `top down' approach but meanwhile the difference is that  we
are explicitly building up the different layers of
infrastructure of the noncommutative geometry.

\section{Local formulae in tensor product gauge}

We are now going to give how the quantum Riemannian
geometry  above looks in the `tensor product gauge'  but only in outline --
for details see \cite{Ma:rie}  and  explicitly \cite{Ma:rieq}  with general
calculi. Also if you do not like Hopf algebras we  give an even more explicit
formulation in the next section.
Thus, let $M$ be an algebra and $H$ a Hopf algebra. We take
$P=M\tens H$, which is the `tensor product gauge' for a trivial
principal bundle. We have differential calculi on $\Omega^1(H)$,
$\Omega^1(M)$ defined by subbimodules $N_M, N_H$ say. We take for
$\Omega^1(P)$ the calculus defined by sub-bimodule
\eqn{NP}{
N_P=N_M\tens H\tens H+M\tens M\tens N_H+\Omega^1 M\tens \Omega^1H}
one has to work to show that we fulfill the conditions for a
quantum principal bundle and that in fact
$P\Omega^1(M)P=\Omega^1(M)P=P\Omega^1(M)=\Omega^1(M)\tens H$
which then allows all formulae in the global theory to be lowered
to working just with $M$ and $H$. For the global picture of
curvature and torsion etc., one also has to specify the global
$\Omega^2(P)$ from $\Omega^2(M)$ and the $\Omega^2(H)$, etc. but
this can be done for any reasonable  $\Omega^2(M)$ (such as the
bicovariant one if $M$ is itself a quantum group).

We then define a gauge field as any map \eqn{Aloc}{ A:\Omega_0\to
\Omega^1(M)} and can show that one can build up a global
connection on the bundle $P$ from this data. We use
$\Ad:\Omega_0\to \Omega_0\tens H$ to do this. It is very
important since we do not consider global gauge transformations
that one does indeed get a global bundle and connection in this
way\cite{Ma:rieq}. Likewise other data above have local analogues
as follows.

A framing amounts to a $V$-bein, i.e. a comodule $V$ (so $H$
(co)acts on it) and a linear map $e:V\to \Omega^1(M)$ such that
it induces an isomorphism $\Omega^1(M)\isom M\tens V$. If you
choose a basis $\{e_i\}$ of $V$ then what we require is that
every $\alpha\in\Omega^1(M)$ has the form $\alpha=\alpha^ie(e_i)$ for unique
functions $\alpha^i$. A coframing is similarly provided by a $V$-cobein $e^*$
such that  $\Omega^1(M)\isom V^*\tens M$. The vanishing of torsion and
cotorsion correspond to
\eqn{torloc}{\bar D\wedge e=0,\quad D\wedge e^*=0.}
Explicit formulae for $D$ on forms are part of the standard gauge bundle
theory\cite{BrzMa:gau} and look just as in \cite{Ma:rieq}. We
project such equations down to $\Omega^2(M)$. The metric is
$g=e^*\tens_M e$ evaluated on the canonical element of $V^*\tens
V$.

Given a framing and $A$, the covariant derivative is
\eqn{nablaloc}{ \nabla\alpha=\extd \alpha^i\tens_M e(e_i)-\alpha^i
A(\tilde\pi_{\Omega_0} e_i\bz)\tens_M e(e_i\bt),}
where the $\bz$ and $\bo$ refer to the pieces in $H\tens V$ resulting from
applying the coaction to an element of $V$. The
$\tilde\pi_{\Omega_0}$  projects $H$ down to $\Omega_0$. The curvature
corresponds to  \eqn{curvloc}{ F(v)=\extd
A(\tilde\pi_{\Omega_0}v)+A(\tilde\pi_{\Omega_0}v\o)\wedge
A(\tilde\pi_{\Omega_0}v\t),} where the $\o$ and $\t$ refer to the pieces
resulting from the coproduct of $H$. We require for a well-defined  map
$F:\Omega_0\to\Omega^2(M)$ a regularity condition on $A$  which in this gauge
appears as \eqn{regloc}{  A(\tilde\pi_{\Omega_0}q\o)\wedge
A(\tilde\pi_{\Omega_0}q\t)=0,\quad \forall q\in Q_H,} where  $Q_H\subset H$
defines $\Omega_0$ for the bicovariant  differential calculus on $H$. It is
exactly this condition that  drops out with the universal calculus on $H$,
when $Q_H$ is zero. From the curvature we build the Riemann curvature  $R$,
and given a bimodule inclusion $i:\Omega^2(M)\to
\Omega^1(M)\tens_M\Omega^1(M)$ the Ricci tensor,
\eqn{ricciloc}{ {\rm
Ricci}=\<i(R)e(e_i),f^i\>=i(F_A(\tilde\pi_{\Omega_0}
e_i\bz)^{ij}e(e_j)\tens_M e(e_i\bo),} where $F_A=F_A^{ij}e(e_i)\tens_M e(e_j)$
defines its components.

Finally, a Dirac operator is defined by a linear map $\gamma:V\to
\End(W)$ for some other $H$-comodule $W$. It takes the form on
spinors $\psi=\psi\uo\tens\psi\ut\in M\tens W$,
\eqn{diracloc}{ \Dsl\psi= (\del^i\tens \gamma_i)\psi -\psi\uo
A^i(\tilde\pi_{\Omega_0}\psi\ut\bz)\tens
\gamma_i(\psi\ut\bo),} where $\gamma_i=\gamma(e_i)$ and $\extd
 m=(\del^i m)e(e_i)$ defines the
partial derivatives associated to  the $V$-bein, and  $\bz$, $\bo$ are the
pieces in $H\tens W$ of the coaction on $W$.

These quantum group formulae may look unpalatable to anyone not
happy with Hopf algebras. One may also rewrite the coaction of
$H$ as an action of a dual Hopf algebra $U$ of `enveloping
algebra' type which will make them look more familiar. Moreover,
see the next section.

\section{Gravity on finite sets}

We now specialise further to concrete local formulae with lots of
indices, but for simplicity only in the simplest case of finite
sets and framing by a classical group. It is very important
conceptually and also (not so important) historically that this
theory on finite sets is {\em not} something new but merely an
elaboration of a more general theory essentially introduced in
1997 and that can also be specialised in other limits, e.g.
noncommutative or classical geometry.

At this level, let $\Sigma$ be a finite set, $G$ a finite group,
$H=\C[G]$ and $M=\C[\Sigma]$ spanned by delta-functions
$\{\delta_x\}$ for $x\in\Sigma$. It is trivial to see (and
well-known) that a general differential calculus $\Omega^1(M)$
corresponds to a subset \eqn{1formset}{ E\subseteq
\Sigma\times\Sigma-{\rm diagonal},\quad
\Omega^1(M)=\{\delta_x\tens\delta_y|\ (x,y)\in E\}} where we set
to zero delta-functions corresponding to the complement of $E$
and identify the remainder with their lifts as shown. If $f=\sum f_x\delta_x$
is a function with components $f_x$, then $\extd f$ has components $(\extd
f)_{x,y}=f_y-f_x$ for $(x,y)\in E$.

Then a $V$-bein is a vector space on which $G$ acts by $\rho_V$, say, and
a collection of 1-forms
\eqn{Efin}{E_i=\sum_{(x,y)\in  E}E_{i,x,y}\delta_x\tens\delta_y} for each
element of a basis  $\{e_i\}$ of $V$ such that the matrices $\{E_{i,x,y}\}$
are invertible for each $x\in \Sigma$ held fixed. A necessary (and sufficient)
condition for the existence of a $V$-bein is clearly that $E$ is fibered over
$\Sigma$, i.e. for each $x\in\Sigma$ the set $F_x=\{y|\ (x,y)\in E\}$ has the
same size, namely the dimension of $V$. In particular it implies that
the latter is $|E|/|\Sigma|$. A natural `local' class of $V$-beins is
just given by any collection of bijections
\eqn{locbeinfin}{s_x:\{i\}\isom F_x,\quad E_{i,x,y}=\delta_{s_x(i),y},}
but we are not limited to such a class. Similarly a $V$-cobein is a
collection of  1-forms with components $E^{*i}_{x,y}$ with respect to a dual
basis $\{f^i\}$ and with the matrices  $\{E^{*i}_{x,y}\}$ invertible for each
$y\in\Sigma$ held fixed. The metric is \eqn{gfin}{  g=\sum_{(x,y,z)\in
F}g_{x,y,z}\delta_x\tens\delta_y\tens\delta_z,\quad
g_{x,y,z}=E^{*i}_{x,y}E_{i,y,z},} where $F$ is the subset of  $(x,y,z)\in
\Sigma\times\Sigma\times\Sigma$ with $(x,y)\in E$ and  $(y,z)\in E$, i.e. $F$
labels the basis of $\Omega^1(M)\tens_M\Omega^1(M)$.

As immediate from Woronowicz'
paper \cite{Wor:dif}, bicovariant  calculi on $H=\C[G]$ are classified by
nontrivial conjugacy  classes or more generally (reducible ones) by
$\Ad$-stable  subspaces $\CC\subset G$ excluding the group identity $e\in
G$. We denote its elements by $a,b,c$ etc., and identify the quotient with the
corresponding lift so that  \eqn{Omega0G}{ \Omega_0=\{\delta_a|\ a\in\CC\}.}
Then a connection or gauge field with values in the dual  of $\Omega_0$ is a
collection of 1-forms with components  $A_{a,x,y}$. In our case $G$ acts on
$V$ so that it plays the  role of frame transformations in the frame bundle
approach. In  that case $A$ induces a covariant derivative on 1-forms
\eqn{nablafin}{ (\nabla\alpha)_{x,y,z}=
(\alpha^i_y-\alpha^i_x)E_{i,y,z}-\alpha^i_xA_{a,x,y}E_{j,y,z}
\tau^a{}^j{}_i,\quad  \tau^a=\rho_V(a^{-1}-e),} where $\alpha=\alpha^iE_i$
defines the  component functions $\alpha^i$ of a 1-form $\alpha$ in the
$V$-bein basis. The $\tau^a$ are the matrices for the action of  the
`braided-Lie algebra' dual of $\Omega_0$ as a subspace of the linear span $\C
G$.

 Next we specify $\Omega^2(M)$ by a bimodule surjection
$\wedge:\Omega^1(M)\tens_M\Omega^1(M)\to \Omega^2(M)$. We assume
that this is done in a manner compatible with $\Omega^2(H)$ so as
to fit together globally in the bundle. The surjection on
any $f\in\Omega^1(M)\tens\Omega^1(M)$ with components $f_{x,y,z}$ as above, is
necessarily of the form
\eqn{wedgefin}{(\wedge
f)_{x,\alpha,z}=\sum_{y\in F_{x,z}}f_{x,y,z}p^y{}_\alpha,\quad
F_{x,z}=\{y\in\Sigma|\ (x,y,z)\in F\}}
for a family of surjections $p:\C F_{x,z}\to V_{x,z}$ to some vector spaces
$V_{x,z}$, obeying $p(1,1,\cdots,1)=0$ when $(x,z)\notin E$ with $x\ne z$ (this
is so that $\extd^2=0$). Chosing a basis $\{e_\alpha\}$ for the latter, we
write the $p$ explicitly as a family of rectangular matrices with each row
summing to 0 in the stated case.  When $\alpha_{x,y},\beta_{x,y}$ are the
components of 1-forms as above then
\eqn{extd2fin}{ (\extd\alpha)_{x,\alpha,z}=\sum_{y\in
F_{x,z}} (\alpha_{x,y}+\alpha_{y,z}-\alpha_{x,z})p^y{}_\alpha,\quad
(\alpha\wedge\beta)_{x,\alpha,z}=\sum_{y\in
F_{x,z}}\alpha_{x,y}\beta_{y,z}p^y{}_\alpha.}

With such an explicit description of $\Omega^2(M)$, a connection $A$ is
regular in the tensor product gauge if
\eqn{regfin}{  \sum_{ab=q,y}A_{a,x,y}A_{b,y,z}p^y{}_\alpha=0,\quad \forall
q\notin\CC\cup\{e\}.}   Its  curvature is
\eqn{Ffin}{ F_{a,x,\alpha,z}=(\extd
A_a)_{x,\alpha,z}+\sum_{cd=a,y}A_{c,x,y}A_{d,y,z}p^y{}_\alpha-
\sum_{b,y}(A_{b,x,y}A_{a,y,z}+A_{a,x,y}A_{b,y,z})p^y{}_\alpha.}
The actual Riemann tensor if one
wants it is the 2-form valued operator on 1-forms,
\eqn{riemfin}{ R_{x,\alpha,z}{}^i{}_j=F_{a,x,\alpha,z}\tau^a{}^i{}_j,\quad
R\alpha=\alpha^i R^j{}_i\tens_M E_j.}
Meanwhile, the zero torsion and  zero cotorsion equations are vanishing of
\eqn{torfin}{(\bar D \wedge e)_{i,x,\alpha,z}=(\extd
E_i)_{x,\alpha,z}+\sum_{a,j,y}A_{a,x,y}E_{j,y,z}p^y{}_\alpha\tau^a{}^j{}_i,}
\eqn{cotorfin}{(D\wedge e^*)^i_{x,\alpha,z}=(\extd
E^{*i})_{x,\alpha,z}+\sum_{a,j,y}E^{*j}_{x,y}A_{a,y,z}p^y{}_\alpha\tau^a{}^i{}_j.}

Also, given a `lift' $i$, which means a collection of inclusions $i:V_{x,z}\to
\C F_{x,z}$ or rectangular matrices $i^\alpha{}_y$, preferably such that
$p\circ i=\id$ (in which case $i\circ\wedge$ is a projection operator splitting
$\Omega^1(M)\tens_M\Omega^1(M)$ into something isomorphic to $\Omega^2(M)$
plus a complement), we have an interior product and, in particular, a Ricci
tensor
\eqn{riccifin}{ {\rm
Ricci}_{x,y,z}=i(F_a)_x^{ij}E_{j,x,y}E_{k,y,z}\tau^a{}^k{}_i.} Here
$i(F_a)_{x,w,z}$ in $\Omega^1(M)\tens_M\Omega^1(M)$ is as in (\ref{Ffin}) but
with $\pi^y{}_w=p^y{}_\alpha i^\alpha{}_w$ in place of $p^y{}_\alpha$
written there. These are projections if $i$ is a proper lift, but this is
not strictly required. We then convert to $V$-bein components $i(F_a)^{ij}$ as
usual. One can write this more explicitly in terms of $A,e$.

Finally, gamma-matrices are a collection of matrices  $\gamma_i$
acting on spinors $\psi$ which are functions with values in a vector space $W$
on which $G$ acts by $\rho_W$, say. The associated Dirac  operator is
\eqn{dirfin}{ \Dsl=\del^i\gamma_i-A_a^i\gamma_i \tau^a_W,\quad
\tau_W^a=\rho_W(a^{-1}-e).}

This is the elementary elaboration for finite sets
of the general  theory in preceding sections. Of course, a simple case is
where  $\Sigma=G$ and we use the same bicovariant calculus on both.  We take
$V=\Omega_0$  itself with the coadjoint action of $G$, so that
the indices $i,j,k$ run over the same range $\CC$ as the $a,b,c$. As
proposed in \cite{Ma:rie} we take the $V$-bein to be the quantum group
Maurer-Cartan form $e$ from \cite{Wor:dif}. There is also a natural
braided-Killing form $\eta$ associated in \cite{Ma:lie} to $\Omega_0^*$ as a
braided-Lie algebra, which we use to define $e^*=e\circ\eta$. The subset $E$
above, the exterior derivative, the components of Maurer-Cartan form and the
matrices $\tau^a$ are explicitly,
\eqn{Gdif}{E=\{(x,y)|\
x^{-1}y\in\CC\},\quad \extd f=(\del^if)E_i,\quad \del^i=R_i-\id,\quad
E_{i,x,y}=\delta_{xi,y},\quad
\tau^a{}^i{}_j=\delta^i{}_{a^{-1}ja}-\delta^i{}_j,}
where $R_i$ is right-translation by $i$. Most of this is the standard starting
point for any work on noncommutative differential geometry of quantum groups,
translated into our above notations.  Also, note that in all computations one
can either work explicitly with matrices $\alpha_{x,y}$ etc. which is like
`spacetime coordinates' or work more algebraically with $V$-bein components
$\alpha^i_x$ and the abstract relations in $\Omega^1(M),\Omega^2(M)$ (or a
mixture of the two). The  conversion in the group case is particularly easy
because of the translation-invariant form of $E_i$.

\begin{theorem}\cite{Ma:rieq} For
$G=\Sigma=S_3$, the permutation group on 3  elements, and $\CC$ the maximal
(order 3) conjugacy class (i) the  braided Killing form is
$\eta^{ij}=3\delta^{ij}$ and defines the metric
$g_{x,y,z}=3\delta_{x^{-1}y,y^{-1}z}$, (ii) there is a  unique torsion-free
and cotorsion-free regular or
`generalised Levi-Civita' connection for this
metric, given by $A_{a,x,y}=\delta_{xa,y}-{1\over 3}$. \end{theorem}

We also look in \cite{Ma:rieq} at two natural `lifts' in the finite group case.
One is the Woronowicz lift whereby $i=\id-\Psi$. The kernel of this map is
precisely the elements in $\Omega^1(M)\tens_M\Omega^1(M)$ set to zero by
$\wedge$, so we can identify $\Omega^2(M)$ with its image. Since $\Psi^2\ne
\id$ for a nonAbelian $G$, it is not a precise lift in the sense of
$i\circ\wedge$ a projection. Another is an actual lift constructed in
\cite{Ma:rieq} giving a projection operator for all finite groups. In either
case one finds\cite{Ma:rieq} \eqn{ricciS3}{ {\rm Ricci}_{x,y,z}=-\mu
(g_{x,y,z}-1),} where $\mu$ is a positive constant. The further contraction,
the scalar curvature, is also constant and negative. Note that this is like a
sphere in our preferred conventions for Ricci (where we contracted what would
classically be the first and third indices of $R$, i.e. minus the usual
conventions).

Finally, there are natural $2\times 2$ equivariant $\gamma$-matrices obeying
\eqn{gammaS3}{ \{\gamma_i,\gamma_j\}+{2\over  3}(\gamma_i+\gamma_j)={1\over
3}(\delta_{ij}-1)} and the  corresponding Dirac operator turns
out to be\cite{Ma:rieq}  \eqn{diracS3}{\Dsl=\del^i\gamma_i-1.} This is just the
beginning  (the natural Riemannian structure on $S_3$). Clearly we have the
machinery above to consider the moduli of all $(A,e,e^*)$ and  solve field
equations or Legendre transform action functionals with respect to them,  i.e.
classical and quantum gravity. If one wants a more conventional moduli space
for the metric then one can fix the relationship $e^*=e\circ \eta$ as above
and only vary $(A,e)$. We will consider the moduli  space elsewhere.

\section{Quantum measurement and algebra
bundles}
As mentioned in the introduction, there is one further generalisation beyond
the case of quantum group (or classical group) framing  that is needed for a
fully comprehensive theory of Riemannian geometry that would apply in principle
to real-world quantum systems. That is to the level of only coalgebras in
place of $H$ above, or to algebras in place of its `enveloping type' dual
Hopf algebra $A$. Such a theory more or less exists as the theory of
(co)algebra bundles in \cite{BrzMa:coa} at the gauge theory level and
\cite{BrzMa:geo} at the frame bundle level. The basic idea is the following.
Given a classical or quantum bundle as above, one can make a cross product
\eqn{crossalg}{ X=P\lcross A}
by the action of $A$ on $P$. Like all cross products this factorises into $P$
and $A$ as subalgebras. The generalisation is to consider more general
algebras $X$ that factorise $X=PA$ as vector spaces, where $P,A$ are
subalgebras. In addition we require a map $e:A\to P$ with certain
properties corresponding to the exact sequence that defined a
global bundle in Section~3. It turns out then that most of the global gauge
theory and Riemannian geometry above, i.e. connections as projections of
$\Omega^1(P)$, associated bundles, frames, etc. then goes through even though
$A$ is not a Hopf algebra at all.

To give some idea of the constructions in \cite{BrzMa:coa}\cite{BrzMa:geo},
first note that when $X$ factorises there is an induced `reordering map'
$\Psi:A\tens P\to P\tens A$ called the entwining or factorisation structure,
and when $e$ is given obeying   \eqn{eab}{e(ab)=\cdot(\id\tens e)\Psi(a\tens
e(b)),\quad e(1)=1,} one has an action
\eqn{eact}{ a\la p=\cdot(\id\tens e)\Psi(a\tens p)}
of $A$ on $P$. This is such that $P$ acting on itself and this action of $A$
fit together to an action of $X$ on $P$. Next, we define $M$ as the `fixed
subalgebra' by
\eqn{Mmeas}{ M=\{m\in P|\ a\la m=e(a)m,\quad \forall a\in A\}
=\{m\in P|\ a\la(pm)=(a\la p)m,\quad \forall a\in A,\ p\in P\}.}
We can then view the action of $A$ more precisely as the coaction of a
coalgebra $C$ dual to $A$ and define ${\rm ver}:\Omega^1P\to P\tens
\ker\eps$ as in the exact sequence before, where $\ker\eps$ are the elements
of $C$ that vanish on $1$. This is for the universal differential calculus. The
theory with general differential calculus is not fully elaborated but exists in
principle; we need a quotient $\Omega_0$ of $\ker\eps$ with suitable
properties and the similar exact sequence. Not that these conditions relate to
the differential geometric structure of the gauge group and are therefore
important for the global bundle structure but not at the local level. For
example, for finite sets above we needed for most of the theory only that
$\Omega_0$ was some chosen quotient and none of its properties (until we took
$G$ also as spacetime). Also Ricci and the Dirac operator are not yet
elaborated in this context, but are not a problem in principle.

The first point of interest for conventional quantum theory is that the
converse is also true: if $A$ is an algebra acting on the vector space of some
algebra $P$, define $M$ by the second version in (\ref{Mmeas}),
let $e(a)=a\la 1$ and suppose that
we have the exact sequence or `local triviality' condition. Then one gets a
bundle and a map $\Psi$ leading to a generalised braided tensor product
algebra \eqn{Xpsi}{ X=P{\underline\tens}_\Psi A}
factorising into $P,A$ and still represented on $P$ (with $P$ acting by
multiplication). Here the product $(q\tens a)(p\tens b)$ is to use $\Psi$ to
take $a$ past $p$ and then multiply in $A$ and in $P$.

We now outline how this mathematics could form the basis of a quantum theory of
measurement. Usually we have some fairly naive `postulates' that a
measurement is deemed to have taken place when a wave-function is projected
into an eigenstate of some chosen operators. The problem is that from a larger
point of view the measurement itself should be a process in a larger quantum
system and only appearing like the above in some idealisation. In the larger
system we would need to identify the algebraic structures which would become
the macroscopic parts of the measuring apparatus and the remainder the
quantum part in that idealisation. Such identifications or correspondences is
exactly the task of quantum geometry. Now in quantum mechanics  this
geometrical structure is usually buried in the Hilbert space, which is
typically of the form $L^2(C)$ for some `configuration space' $C$. It is not
just some abstract Hilbert space (these are all isomorphic) because the
structure of $C$ goes into the construction of the Hamiltonian and other
correspondences between the classical and quantum theory. Forgetting about the
$L^2$ completion, we are therefore in the situation where our quantum system
is an algebra $A$ acting on the vector space of some other {\em algebra} $P$.
The algebraic structure of $P$ carries the geometric
structure of $C$. This exactly the situation above, so we have as an immediate
application of \cite{BrzMa:coa}\cite{BrzMa:geo},

\begin{corol} (i) For any quantum system $A$ acting on (the Hilbert
space completion of) an algebra of wave-functions $P$, we have a subalgebra
$M\subseteq P$ defined by (\ref{Mmeas}), the `superselection algebra' of the
quantum system. We say that the system is reduced if $M$ is trivial. (ii) If
the exact sequence condition holds so that we have a bundle (we say that the
quantum system is {\em Galois}) then we have an extended quantum system
$P{\underline\tens}_\Psi A$ containing both the quantum system $A$ and the
algebra $P$ of functions on the configuration space, represented together on
$P$. \end{corol}

In the Galois case the bigger algebra contains both what was to be measured
(the quantum system) and the classical system $P$ of functions on the
configuration space wherin measurements lie (in the
sense that the norm of the value at $c$ is the probability density to measure
$c\in C$). As well as a method of constructing new quantum systems from old
(on the same Hilbert space), this is also a method of constructing quantum
systems in the first place from a representation of one part $A$ that we would
like to include in another part $P$; even if both parts are classical the
result $X$ can be quantum. Actually, this is a very typical construction for
quantum systems, where we typically seek to represent momentum modes on the
configuration space and then throw in the configuration space coordinates
acting by multiplication as well (the extended system). The theory above
gives a mechanism for such a construction based just on the act of representing
one part $A$ in another part $P$ subject to a nondegeneracy (Galois) condition.
A concrete example is a
dynamical system where the algebra $A$ generated by a group $G$ acts on the
algebra of functions $P$ on a space $C$. The space $C$ foliates into orbits and
$M$ is the algebra of functions on the space $C/G$ of orbits. The system is
reduced if there is just one orbit. Otherwise $M$ describes the superselection
sectors of the theory. Here $X=P\lcross A$ is what is usually called the
quantisation of the dynamical system, represented on $P$. The geometric
picture here is also clear; $C$ in nice cases is a bundle over $C/G$.
This shows how the above theory is already a useful model of {\em
conventional} quantum mechanics. We can ask which familiar systems are
reduced? Which are Galois? Moreover, in the latter case we have in principle
all the machinery of bundles, etc. to apply quite generally, with $M$ in the
role of base manifold even for quantum systems not obviously of geometric
origin.

On the other hand there is no reason why $P$ (or $A$ or $M$) should be
commutative, i.e. all concepts in the formulation above work when the system
is already noncommutative or quantum. We therefore have the possibility of
quantising a system with `wave-functions' actually operators in some already
quantised system $P$.  Thus we could ask if one part $A$ of a quantum system is
represented on another part  $P$ and build the extended system. This is a
necessary (though not sufficient) part of the solution of the measurement
problem because it begins to formulate procedures in a manner applicable
before making any idealisations to the special form of a classical
configuration space, i.e. measuring quantum theory within quantum theory or an
intrinsic theory of measurement involving some parts `measuring' other parts.
Moreover, we still have the geometric picture but now as quantum bundles,
quantum Riemannian geometry etc. for such systems.  Finally, we also have the
first part of the theory above which goes in the other direction: a
theory of factorisations of a bigger system $X$ into a part that we would like
to be a smaller quantum system $A$ and a part that we would like to be (only
approximately if noncommutative) a classical or macroscopic system $P$.

Putting these ideas together we arrive at a relative theory of quantum
measurement which we propose along the following lines.

\begin{enumerate}
\item  A quantum system $A$ means a $*$-algebra.
\item A relative `wave function' state is an element $p$ of another possibly
quantum system $P$ on which $A$ is represented as a vector space. We call
$1\in P$ the `character state'.
\item We call $\<a\>_p=p^*(a\la p)$ the
`relative expectation value of $a$ in state $p$'. More generally, a relative
expectation is any positive map $\phi:A\to P$.
 \item If $P$ is itself
represented in some other system $Q$ then $\<\<a\>_p\>_q$ is the expectation
of the above expectation in state $q$, etc.
\item When the representation of
$A$ in $P$ is Galois, we may view them as subalgebras of a bigger quantum
system $X=PA$.
\item Conversely, if a quantum system $X$
factorises as $X=PA$ and a relative expectation $e:A\to P$ obeying (\ref{eab})
is specified then $A$ is represented in $P$ and $e(a)=\<a\>_1$.
\end{enumerate}

This relative theory unifies both the Schroedinger and the $C^*$-algebra
type attitudes towards quantum mechanics. Thus, in one extreme, take
$P=\C$. Then a relative expectation is just a positive linear map $\phi:A\to
\C$, i.e. states from a $C^*$-algebra type point of view. A representation
in $\C$ is determined by a character $e:A\to\ C$ and $\<a\>_1=e(a)$. We view
$\C$ as the `end of the line' in that we do not expect to represent this
elsewhere, and can speak (glibly) about actual expectations in this case.
However, our view in general is a Bayesian one in which all probability is
relative. In another special case take $P$ functions on a classical space $C$.
Then wave functions are usual wave functions. $\<a\>_p$ is also a
wave-function, the expectation density. A further representation of the
classical system $P$ in $Q=\C$ is a point $c\in C$ and $\<p\>_1=p(c)$. More
generally, an  expectation state $\phi:P\to \C$ is a convex combination of
points in $C$ or a probability measure, more precisely. Then composing
$\<a\>_p$ with this gives the actual expectation value of $a$ in the usual
sense. For example, if we chose $\phi$ to be the uniform probability
distribution then $\phi(\<a\>_p)$ is the usual quantum mechanical expectation.
It corresponds to the `all other things being equal' assumption in the final
step to get from a relative probability to a so-called absolute one. We would
also expect a relative GNS-type construction based on quotienting the regular
representation of $A$ relative to a general state $A\to P$, etc. This should
go hand in hand with some notion of relative completions needed to make the
above precise but without assuming that every algebra is normed (or a $C^*$
algebra) to begin with (this is optional but would be in keeping with the
relative philosophy).

While not a full theory of quantum measurement, we have indicated here the
start of a general framework. Moreover, this framework has all the tools of
noncommutative Riemannian geometry, gravity etc. as explained above, at our
disposal. Putting these together one would expect a full resolution of the
link between measurement, entropy and gravity. Also note that, while
not Galois, important algebra factorisations arose in the theory of Hopf
algebra factorisations. In those models one has dual Hopf algebras related to
Hopf algebra duality and observable state duality\cite{Ma:pla}. This too should
be generalised and resolved using the above formulation.

\subsection*{Acknowledgements} It is a pleasure to thank the
organisers of the Euroconference for a thoroughly enjoyable and
stimulating event.

\end{document}